\documentclass[final]{elsarticle}

\usepackage[left= 3cm, right=3cm]{geometry}
\usepackage{natbib}

\usepackage{pstricks, pst-plot}	
\usepackage{graphicx} 
\usepackage{wrapfig}  
\usepackage[figuresright]{rotating}

\usepackage{amsmath}
\usepackage{amssymb}
\usepackage{hyperref}

\usepackage[english]{babel}
\usepackage{blindtext}
\usepackage{caption}
\usepackage{subcaption}
\usepackage{algorithm}
\usepackage{algpseudocode}
\usepackage{float}
\newcommand{\sign}{\text{sign}}

\begin{document}
\begin{frontmatter}%

\title{Reinforcement Learning for Batch Bioprocess Optimization}
 \author[a,b]{P. Petsagkourakis}
 \author[c]{I.O. Sandoval}
 \author[d]{E. Bradford}
 \author[a,e]{D. Zhang\corref{cor1}}
 \ead{dongda.zhang@manchester.ac.uk}
 \cortext[cor1]{Corresponding author}
 \author[e]{E.A. del Rio-Chanona\corref{cor2}}
 \ead{a.del-rio-chanona@imperial.ac.uk}
 \cortext[cor2]{Corresponding author}

 \address[a]{School of Chemical Engineering and Analytical Science,The University of Manchester, M13 9PL, UK}
 \address[b]{Centre for Process Systems Engineering (CPSE), Department of Chemical Engineering, University College London, Torrington Place, London WC1E 7JE, United Kingdom}
 \address[c]{Instituto de Ciencias Nucleares, Universidad Nacional Aut\'{o}noma de M\'{e}xico, A.P. 70-543, C.P. 04510 Ciudad de M\'{e}xico, Mexico}
 \address[d]{Department of Engineering Cybernetics, Norwegian University of Science and Technology, Trondheim, Norway}
 \address[e]{Centre for Process Systems Engineering (CPSE), Department of Chemical Engineering, Imperial College London, UK}
%
%
%
%
\begin{abstract}
Bioprocesses have received a lot of attention to produce clean and sustainable alternatives to fossil-based materials. However, they are generally difficult to optimize due to their unsteady-state operation modes and stochastic behaviours. Furthermore, biological systems are highly complex, therefore plant-model mismatch is often present. To address the aforementioned challenges we propose a Reinforcement learning based optimization strategy for batch processes.
\par

In this work we applied the Policy Gradient method from batch-to-batch to update a control policy parametrized by a recurrent neural network. We assume that a preliminary process model is available, which is exploited to obtain a preliminary optimal control policy. Subsequently, this policy is updated based on measurements from the \emph{true} plant. The capabilities of our proposed approach were tested on three case studies (one of which is nonsmooth) using a more complex process model for the \emph{true} system embedded with adequate process disturbance. Lastly, we discussed advantages and disadvantages of this strategy compared against current existing approaches such as nonlinear model predictive control.
\end{abstract}
\end{frontmatter}
%
%
%
%
\section{Introduction}
There has been a global interest in using sustainable bio-production systems to produce a broad range of chemicals to substitute fossil derived synthetic routes \citep{Brennan2010,Harun2018}. Bioprocesses exploit microorganisms to synthesize platform chemicals and high-value products by utilizing different means of resources \citep{Jing2017}. Compared to a traditional chemical process, a biochemical process is highly complex due to the intricate relationships between metabolic reaction networks and culture fluid dynamics \citep{Rio-Chanona2018}. As a result, it is difficult to construct accurate physics-based models to simulate general large-scale biosystems, and plant-model mismatch is inevitable. Furthermore, the behaviour of bioprocesses is often stochastic in the macro-scale, because biological metabolic pathways are sensitive even to mild changes of process operating conditions \citep{Zhang2015,Thierie2004}. Therefore, developing control and optimization strategies for bioprocesses remains an open challenge. Given these critical limitations in mechanistic models, in this work we propose a data-driven approach to address this issue.

We must seek a strategy that can optimize our process and handle both the system's stochastic behaviours (e.g. process disturbances) and plant-model mismatches. It is here that we have opted to use {\it Reinforcement learning} and more specifically, {\it Policy Gradients}, as an alternative to current existing methods. Reinforcement learning (RL) was developed to address nonlinear and stochastic optimal control problems \citep{Bertsekas:2000:DPO:517430}. Two main branches have been established on how to solve dynamic optimization problems via RL.

The first branch is based on {\it Dynamic Programming} (DP), hence termed Approximate Dynamic Programming (ADP). DP relies on the Hamilton-Jacobi-Bellman equation (HJBE), the solution of which becomes intractable for small size problems with nonlinear dynamics and continuous state and control actions. Because of this, past research has relied on using ADP techniques to find approximate solutions to these problems ~\citep{Sutton}.

The second branch is to use {\it Policy Gradients}, which directly obtains a policy by maximizing a desired performance index. This approach is particularly well suited to deal with problems where both the state space and the control space are continuous. Given the advantages that Policy Gradients can offer when confronted with bioprocess optimization, we have adopted this approach in the current work. Policy Gradient methods, along with their benefits, are further explained in Section ~\ref{PGM}. 

\subsection{Related Work}
Given that chemical engineers have always dealt with complex and uncertain systems there have been several approaches that address specific instances of the aforementioned problems, we highlight some related previous work in the sections below.

One example to track stochastic batch-to-batch systems is Iterative Learning Control (ILC) which was initially introduced for robot manipulators \citep{Arimoto}, and later implemented by the process control community \citep{XU19991535}. ILC deals with the problem of tracking the control performance in batch processes given a reference trajectory for runs that last a fixed time, and where the process state is reset to the same value at the start of each run. An overview of ILC strategies in process control can be found in  \citep{LEE20071306}. 

Real-time optimization (RTO) is another method that deals with uncertain processes. The main idea is to represent the process dynamics by a nonlinear input/output mapping where the disturbances are explicitly accounted for. This mapping is then used to optimize some desired performance index. For the interested reader, further details can be found in \cite{Bonvin:28363} and \cite{Chachuat2009}. A recent review on this topic can be found in \cite{Marchetti2016}. For the dynamic systems, these methodologies are usually referred to as Dynamic Real-time optimization (DRTO) which is closely related to NMPC. More details can be found in \cite{Rawlings,Rossi2019}.

Another technique that deals with stochastic systems is model predictive control (MPC), and its extension to nonlinear systems, NMPC. NMPC has a vast variety of methods that can incorporate uncertainty or maintain properties under the presence of stochastic environments. The most common paradigms are the {\it stochastic} NMPC~\cite{Mesbah2016} and the {\it Robust} NMPC~\cite{Bemporad1999}, where the former incorporates the uncertainty by minimizing (usually) the expectation of the objective function, whilst the latter approach solves a min-max optimization by minimizing the worst case scenario of the uncertainty. Both of these approaches require knowledge regarding the nature of the uncertainty in order to proceed. 

There are different approaches that have been proposed for NMPC frameworks, including scenario~\cite{Bernardini2012} based multi-stage schemes for nonlinear systems~\cite{Lucia2012,Danish2018}, where stochastic programming is utilized and future information is incorporated in an adaptive manner.
Another approach is the use of Gaussian processes~\cite{8550249,Bradford2018} or using (generalized) polynomial chaos expansions~\cite{Kim2013} to model effectively the uncertainties of the process. In the case where no proper information for the uncertainty is available, e.g. there is not enough data to conduct uncertainty quantification, optimal control is explored using the nominal linear or nonlinear available model. 
In terms of solution procedures for the dynamic optimization problem, it is common to use a direct approach after parametrizing and discretizing the control inputs~\cite{Vassiliadisprob21994} or the system dynamics~\cite{biegler2010nonlinear} resulting in a nonlinear programming problem. Although much less common, indirect approaches can also be used, where the necessary conditions of optimality are solved explicitly~\cite{Bonvininderect}. If no information on structural information is known, conservative assumptions can be made in order to establish stability conditions~\cite{Feller2016,Petsagkourakis2018, petsagkourakis2019}.

Reinforcement Learning (in an Approximate Dynamic Programming philosophy), has lately caught significant attention for chemical process control. For example, in \cite{Lee2005} a model-based strategy and a model-free strategy for control of nonlinear processes were proposed, in \cite{Peroni2005} ADP strategies were used to address fed-batch reactor optimization, in \cite{Lee2006} mixed-integer decision problems were addressed with applications to scheduling. In \cite{Tang2018} with the inclusion of distributed optimization techniques, an input-constrained optimal control problem solution technique was presented, among other works (e.g. \cite{Chaffart2018}, \cite{Shah2016}). All these approaches rely on the (approximate) solution of the HJBE, and have been shown to be reliable and robust for several problem instances. 

In this paper, we present another take on RL, that of using Policy Gradients. Policy Gradient methods directly estimate the control policy, without the need of a model, or the solution of the HJBE, its advantages are highlighted in the following section.

In addition to the above, for recent reviews of Machine Learning and Artificial Intelligence applied to chemical engineering the reader is referred to \cite{Lee2018} and \cite{Venkatasubramanian2019}. A shorter review focused on modelling bioprocesses with ML tools can be found in \cite{DelRioChanona2018Review}.

\subsection{Motivation}
The process systems engineering community has been dealing with stochastic batch-to-batch systems for a long time. For example, nonlinear dynamic optimization and particularly NMPC are a powerful methodology to address uncertain dynamic systems, however there are several properties that  make its application less attractive.
All the approaches in NMPC require the knowledge of a detailed model that describe the system dynamics, and stochastic NMPC additionally requires an assumption for the uncertainty quantification/propagation. Furthermore, the online computational time may be a bottleneck for real time applications since a (possibly) nonlinear optimization problem has to be solved. 

In contrast, RL directly accounts for the effect of future uncertainty and its feedback in a proper ‘closed-loop’ manner, whereas conventional NMPC assumes open-loop control actions at future time points in the prediction, which can lead to overly conservative control actions~\cite{Lee2005}.
In addition, policy gradients can establish a policy in a model-free fashion and excel at on-line computational time. This is because the online computations require only evaluation of a policy, since all the computational cost is shifted off-line. 

As mentioned previously, Real-time optimization (RTO) has been used to address many instances of batch-to-batch problems. Interestingly, some recent approaches have suggested a hybrid modeling strategy, where function approximates are used in conjunction with a pre-existing model \cite{del-Rio-ChanonaDycops, Gao2015}. From some perspectives these recent algorithms could be thought of as model-based Reinforcement learning approaches. However, there is not yet a clear consensus on how to address problems with plant-model mismatch, measurement noise, and disturbances in an RTO framework. 

In terms of previous RL approaches in chemical engineering to address process control and optimization, they have relied on action-value methods (e.g. Q-learning, solution of the HJBE). However, to address continuous nonlinear action domains Policy Gradient methods present several advantages: 
\begin{itemize}
    \item In Policy Gradient methods, the approximate policy can naturally approach a deterministic policy, whereas action-value methods (that use epsilon-greedy or Boltzmann functions) select a random control action with some heuristic rule  \cite{Sutton}.
    \item Although it is possible to estimate the objective value of state-action pairs in continuous action spaces by function approximators, this does not help choose a control action. Therefore, on-line optimization over the action space for each time-step should be performed, which can be slow and inefficient. Policy Gradient methods work directly with policies that emit probability distributions, which is much faster and does not require an online optimization step.
    \item  Policy Gradient methods are guaranteed to converge at least to a locally optimal policy even in high dimensional continuous state and action spaces, unlike action-value methods where convergence to local optima is not guaranteed \cite{Sutton}.
    \item Policy Gradient methods enable the selection of control actions with arbitrary probabilities. In such cases, the best approximate policy may be stochastic \cite{Sutton}.
\end{itemize}

Due to the above advantages, in this work we propose an optimization strategy that uses a Policy Gradient algorithm to optimize batch-to-batch bioprocesses. This work extends our proposed methodology in~\cite{PETSAGKOURAKIS2019919}, the new approach presents a much faster adaptation time by implementing {\it transfer learning} for the efficient adaptation of the policies. Additional more complex case studies and a comparison against NMPC are included. Furthermore, we exemplify both approaches (NMPC and our approach) in a system described by a nonsmooth differential equation model. The difficulty for nonsmooth models is highlighted in~\cite{STECHLINSKI201852}.


\section{Methodology}
\subsection{Problem Statement}
In this work, we assume that the system's dynamics are given by an (generally) unknown probability distribution, following a Markov process:
\begin{equation}
    \textbf{x}_{t+1}\sim p(\textbf{x}_{t+1}|\textbf{x}_t,\textbf{u}_t)
\end{equation}
This system can be approximated by the following  discrete time stochastic
nonlinear system represented as a state-space model:
\begin{equation}
\begin{split}
    &\textbf{x}_{t+1} = f(\textbf{x}_{t},\textbf{u}_{t},\textbf{d}_{t})\\
\end{split}
\end{equation}
where $t$ represents the discrete time, $\textbf{x}_t \in \mathbb{R}^{n_x}$ is the vector of states, $\textbf{u}_t \in \mathbb{R}^{n_u}$ is the vector of inputs, $\textbf{d}_t \in \mathbb{R}^{n_d}$ is the the vector of process disturbances, and $f(\cdot)$ are the nonlinear dynamics of the physical system.

Our strategy seeks to find the optimal policy for a batch process under the presence of disturbances and measurement noise. Then, the problem can be written as an Optimal Control Problem (OCP):
\begin{equation}
\mathcal{P}(\pi(\cdot)):=\left \{\begin{aligned}
        &\max_{\pi(\cdot)} \mathbb{E}[J(\textbf{x}^k_{t},\textbf{u}^k_{t})]\\
    &\text{s.t.}\\
    &\textbf{x}^k_0 = \textbf{x}^k(0)\\
    &\textbf{x}^k_{t+1} = f(\textbf{x}^k_{t},\textbf{u}^k_{t},\textbf{d}^k_{t})\quad \forall t \in \left\{1,...,T-1\right\} \\
    &\textbf{u}_t \sim \pi(\textbf{x}_t^k)\\
    &\textbf{u}\in\mathbb{U}\\
    &\text{given}\\
    &\textbf{x}^j_t \quad \forall j \in \left\{0,...,k-1\right\}  \quad \forall t \in \left\{1,...,T\right\}\label{eq:OCP}
\end{aligned}\right.
\end{equation}
The objective is to maximize the expectation of an economic criterion $J$, where $k$ is the current batch, while $j$ refers to previous batch realizations. Additionally, the optimization problem (\ref{eq:OCP}) searches for a set of functions $\pi(\cdot)$ that maps the probability of $u_t^k$ given $x_t^k$. Notice that in problem (\ref{eq:OCP}) we make no assumptions about the nature of $\textbf{d}$. Even in the case where the dynamics are fully known, the solution of problem (\ref{eq:OCP}) may be intractable for medium size systems. 

\begin{figure}[H]
    \centering
    \includegraphics[scale = 0.3]{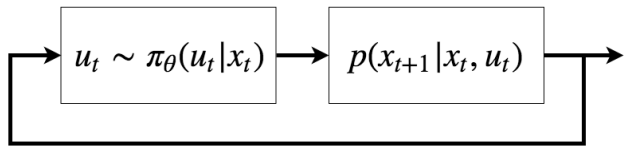}
    \caption{Representation of interaction between policy and the physical system}
    \label{fig:Close-loop}
\end{figure}

To overcome this limitation a novel strategy is proposed, where a policy $\pi_\theta(\cdot)$, parametrized by the parameters $\theta$, is constructed that maximizes the expectation of a performance index $J$. The states at the time $t+1$ are assumed to be given by the probability density $p(\textbf{x}_{t+1}|\textbf{x}_t , \textbf{u}_t)$. The interaction with the policy can be depicted as a closed-loop, see Fig.~\ref{fig:Close-loop}.
Let $\pmb{\tau}$ denote a joint random variable of states and controls defining a trajectory with a time horizon $T$:   
$\pmb{\tau} = (\textbf{x}_0,\textbf{u}_0,R_0,...,\textbf{x}_{T-1},\textbf{u}_{T-1},R_{T-1},\textbf{x}_T,R_T)$, the performance index being 
\begin{equation}
    J(\pmb{\tau}) = \sum_{t=0}^T\gamma^{t}R_t(\textbf{u}_t,\textbf{x}_t)
\end{equation}
where $\gamma\in (0,1]$ is the {\it discount factor} and $R_t$ a given reward at the time instance $t$ for the values of $u_t$, $x_t$.
We represent the probability density of a trajectory as:
\begin{align}
p(\pmb{\tau}|\theta) = \hat{\mu}(\textbf{x}_0) \prod_{t=0}^{T-1} \left[\pi(\textbf{u}_t|\textbf{x}_t,\theta) p(\textbf{x}_{t+1}|\textbf{x}_t,\textbf{u}_t) \right]\label{eq:Eq trajectory}
\end{align}
where $\hat{\mu} (\textbf{x}_0)$ is the probability density of the initial state. We can therefore state the following optimization problem:
\begin{align}
\max_{\pi(\cdot)}\mathbb{E}_{\pmb{\tau} \sim p(\pmb{\tau}|\theta)}[J(\pmb{\tau})]\label{eq:Expectation J}
\end{align}
notice that the process dynamics are implicit in $\pmb{\tau}$. To solve problem (\ref{eq:Expectation J}) we turn our attention to policy gradient methods.
%
%
\subsection{Policy Gradient Methods}\label{PGM} 
Policy gradient methods compute a policy that maximizes the expectation over some objective function ($i.e.$ problem (\ref{eq:Expectation J}) ). They rely on a parametrized policy function $\pi_\theta(\cdot)$ that returns an action $\textbf{u}$ given a state of the system $\textbf{x}$ and a set of intrinsic parameters $\theta$ of the policy. In the case of stochastic policies, the policy function returns the defining parameters of a probability distribution over all possible actions, from which the actions are sampled: 
\begin{align}
\textbf{u} \sim \pi_\theta(\textbf{u} | \textbf{x}) = \pi(\textbf{u} | \textbf{x}, \theta) = p(\textbf{u}_t = \textbf{u}| \textbf{x}_t = \textbf{x}, \theta_t = \theta).
\end{align}
In this work, a Recurrent neural network (RNN) is used as the parametrized policy, which takes (a number of past) states and control actions as inputs and returns the moments of a probability distribution. Then the next control action is drawn from the corresponding probability distribution. For example, if the control actions live in a normal distribution then a mean and a variance are computed, from these mean and variance a control action can be drawn. In this setting, the exploitation-exploration trade-off is represented explicitly by the value of the variance of the underlying distribution of the policy. Deterministic policies can be seen as a limiting case where the variance converges to zero. 

To learn the optimal policy, we seek to maximize our performance metric, and hence we can follow a gradient ascent strategy: 
\begin{align}
\theta_{m+1} = \theta_m + \alpha_m \nabla_\theta \mathbb{E}_{\pmb{\tau} \sim p(\pmb{\tau}|\theta)}[J(\pmb{\tau})] \label{eq:Steepest Ascent}
\end{align}
where $m$ is the current iteration that the parameters are updated (epoch), $\mathbb{E}_{\pmb{\tau}\sim p(\pmb{\tau}|\theta)}[J(\pmb{\tau})]$ is the expectation of $J$ over $\pmb{\tau}$ and $\alpha_m$ is the step size (also termed learning rate in the RL community) for the $m^{th}$ iteration. 
Computing $\hat{J}(\theta) = \mathbb{E}_{\pmb{\tau}\sim p(\pmb{\tau}|\theta)}[J(\pmb{\tau})]$ directly is difficult, therefore we use the \emph{Policy Gradient Theorem} \cite{Sutton:1999:PGM:3009657.3009806}, which shows the following:
\begin{subequations}
\begin{alignat}{3}
\hat{J}(\theta) = \nabla_\theta \mathbb{E}_{\pmb{\tau} \sim p(\pmb{\tau}|\theta)}[J(\pmb{\tau})] &~=~\nabla_\theta \int  p(\pmb{\tau}|\theta)~ J(\pmb{\tau})\text{d}\pmb{\tau}\\\label{eq:expectation}
&~=~\int \nabla_\theta p(\pmb{\tau}|\theta)~ J(\pmb{\tau}) \text{d}\pmb{\tau}\\
&~=~\int p(\pmb{\tau}|\theta) ~ \frac{\nabla_\theta p(\pmb{\tau}|\theta)}{p(\pmb{\tau}|\theta)~} ~ J(\pmb{\tau}) \text{d}\pmb{\tau}\\
&~=~\int p(\pmb{\tau}|\theta) ~ \nabla_\theta \text{log}\left( p(\pmb{\tau}|\theta)\right) J(\pmb{\tau})\text{d}\pmb{\tau}\\
&~=~\mathbb{E}_{\pmb{\tau}} \left[ J(\pmb{\tau}) \nabla_\theta \text{log}\left( p(\pmb{\tau}|\theta)\right) \right] \label{eq:ExpGrad1}
\end{alignat}
\end{subequations}
Notice from (\ref{eq:expectation}) that, $p(\pmb{\tau}|\theta)~ J(\pmb{\tau})$ is an objective function value multiplied by its probability density, therefore, integrating this over all possible values of $\pmb{\tau}$ we obtain the expected value. From there we arrive at (\ref{eq:ExpGrad1}), where, dropping the explicit distribution of $\pmb{\tau}$, gives us an unbiased gradient estimator, (\ref{eq:Steepest Ascent}) now becomes:
\begin{align}
\theta_{m+1} = \theta_m + \alpha_m \mathbb{E}_{\pmb{\tau}}\left[ J(\pmb{\tau}) \nabla_\theta \text{log}\left( p(\pmb{\tau}|\theta)\right) \right] \label{eq:Steepest Ascent 2}
\end{align}
Using the expression for $p(\pmb{\tau}|\theta)$ in (\ref{eq:Eq trajectory})  and taking its logarithm, we obtain:
\begin{align}
\nabla_\theta \text{log}\left( p(\pmb{\tau}|\theta)\right) = \nabla_\theta \sum_{t=0}^{T-1} \text{log}\left( \pi(\textbf{u}_t|\textbf{x}_t,\theta)\right)
\end{align}
Note that since $p(\textbf{x}_{t+1}|\textbf{x}_t,\textbf{u}_t)$ and $\hat{\mu}(\textbf{x}_0)$ are independent of $\theta$ they disappear from the above expression. Then we can rewrite (\ref{eq:ExpGrad1}) for a trajectory as:
\begin{align}\label{eq:PolicyGra}
\nabla_\theta \mathbb{E}_{\pmb{\tau}}[J(\pmb{\tau})] = \mathbb{E}_{\pmb{\tau}} \left[ J(\pmb{\tau}) \nabla_\theta \sum_{t=0}^{T-1} \text{log}\left( \pi(\textbf{u}_t|\textbf{x}_t,\theta)\right)\right]
\end{align}
Notice that expression (\ref{eq:PolicyGra}) does not require the knowledge of of the dynamics of the physical system. However, the above update presents two challenges: the selection of the policy $\pi(\textbf{u}_t|\textbf{x}_t,\theta)$  and the computation of the expectation. To address these possible issues, in this work, recurrent neural networks are used to parametrize the policy of the policy gradient (presented in Section \ref{RNN}), while a Monte-Carlo method is utilized to approximate the expectation (presented in Section \ref{REA}).
%
%



%
\subsection{Recurrent Neural Network}\label{RNN} 

\begin{figure}[H]
    \centering
        \includegraphics[scale = 0.3]{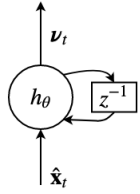}
    \caption{Computational graph of recurrent neural network}
    \label{fig:Fold_RNN}
\end{figure}

Recurrent neural networks, (RNNs)~\citep{Rumelhart1986}, are a type of artificial neural network that has been tailored to address sequential data. RNNs produce an output at each time step and have recursive connections between hidden units. This allows them to have a `memory' of previous data and hence be well suited to model time-series. In general, RNNs can be depicted as a folded computational graph as presented in Fig.~\ref{fig:Fold_RNN}.

\begin{figure}[H]
     \centering
     \includegraphics[scale = 0.3]{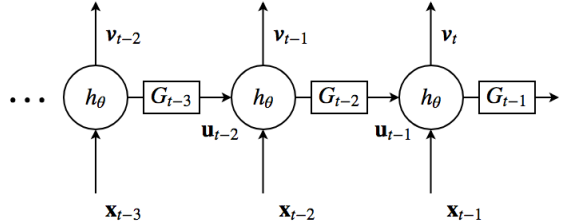}
     \caption{Recurrent neural networks as unfolded computational graph, with one step delay ($z^{-1})$}
     \label{fig:Unfold_RNN}
 \end{figure}

 Fig.~\ref{fig:Fold_RNN} shows how an input $\hat{\textbf{x}}_{t}$ is presented to the network as well as the recursive state of RNN, $\textbf{u}_{t-1}$ and outputs $\pmb{\nu}_t$. A more detailed representation of an RNN is depicted in Fig.~\ref{fig:Unfold_RNN}, which is equivalent to a series of unfolded nodes associated with a particular time instance.
 In Fig.~\ref{fig:Unfold_RNN} we can appreciate that each node receives two inputs $\textbf{x}_{t-i}$ and $\textbf{u}_{t-i-1}$. Generally speaking, the input $\textbf{x}_{t-i}$ corresponds to the data supplied to that node, such as in a traditional artificial neural network (also referred to as feedforward). The unfolded computational graph in this case can be represented as a dynamic system:
 \begin{equation}
     \begin{split}
         \pmb{\nu}_{t} &= h_\theta(\hat{\textbf{x}}_{t},\textbf{u}_{t-1})\\
         \textbf{u}_{t} &= G_{t}(\pmb{\nu}_{t})
     \end{split}
 \end{equation}
where $\hat{\textbf{x}}_t$ is the vector that contains all the external variables for the RNN, $G$ the function that computes the output of the network $u$, and $h_\theta$ represents the layers of the neural networks. Deep structures (which means having more than one hidden layer) can be employed to enhance the performance of the network (deep neural networks) which have been combined with Reinforcement learning recently in~\cite{Mnih2013,Mnih2015}.  Previous realizations of the states $x$ and the controls $u$ are also used as input variables to the network, e.g. $\hat{\textbf{x}}_t = \left[\textbf{x}_{t}^T,\dots,\textbf{x}_{t-N}^T,\textbf{u}_{t-2}^T,\dots, \textbf{u}_{t-N-1}^T \right]^T$, to model dynamic systems. RNNs have previously been applied either as a surrogate model of the process dynamics~\cite{Su1992} or as a parametrization of the agent (control policy)~\cite{Mnih2014}.

In this work, RNNs are applied to parameterize the stochastic policy. We must remark that in theory the Markov decision process does not require RNNs (due to the Markov property), however in practice the use of RNNs can improve the performance of the policy by exploiting additional memory that is provided. In the current work the RNN initially computes the mean and the variance of a multivariate normal distribution where the control actions live. Subsequently, the {\it actual} control action is drawn. Precisely, $\pmb{\nu}_t =\left[ \pmb{\mu}_t,\pmb{\Sigma}_t\right]$ where the $\pmb{\mu}$ and $\pmb{\Sigma}$ are the mean and variance, respectively, and $\textbf{u}_t = G_{t}(\pmb{\nu}_t)$ is substituted by $\textbf{u}_t \sim \mathcal{N}(\pmb{\mu}_t,\pmb{\Sigma}_t)$ making it a stochastic policy. Under the presence of  uncertainty (stochastic in nature) a deterministic policy will fail as the control action will always be the same for the same states since it learns a deterministic mapping from states to control actions at the exact same state. On the contrary, a stochastic policy draws a control action from a probability distribution which can account for stochastic environments.
\begin{figure}[H]
    \centering
        \includegraphics[scale = 0.3]{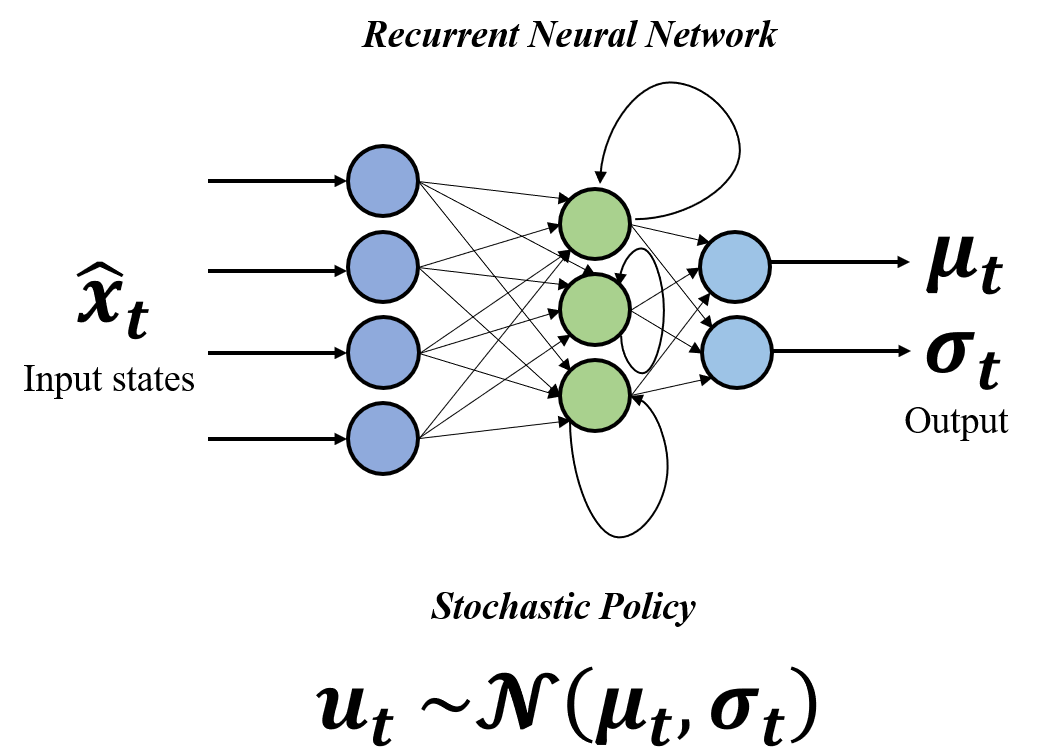}
    \caption{Graphical representation of a stochastic policy network}
    \label{fig:RNN_policy}
\end{figure}
In Figure \ref{fig:RNN_policy}, a schematic representation of a policy network where the stochastic policy follows a Gaussian distribution is depicted. In this figure, we can observe how the states $\hat{x}_t$ are used as input to the network, and how the network computes the mean $\mu_{t+1}$ and standard deviation $\sigma_{t+1}$ for the subsequent time step. Then, the control action $\textbf{u}_{t+1}$ is drawn from the distribution defined by the outputs of the network.
\subsection{Reinforce Algorithm}\label{REA} 
 Given that the parametrized policy used in this work is a RNN, it must be trained to adjust its weights so that the output corresponds to an optimal control action. To this end, we use the steepest ascent strategy mentioned in (\ref{eq:Steepest Ascent}). However, computing the expectation in (\ref{eq:PolicyGra}) can be an intractable problem, and this expression is needed to compute the steepest ascent update ( \ref{eq:Steepest Ascent 2}). Therefore, we propose to use the Reinforce algorithm to compute the policy gradient. The Reinforce algorithm ~\citep{williams1992simple} approximates the gradient of the policy to maximize the expected reward with respect to the parameters $\theta$ without the need of a dynamic model of the process. To compute the expectation we take several sample trajectories and then approximately calculate $\mathbb{E}_{\pmb{\tau}} \left[ J(\pmb{\tau}) \nabla_\theta \sum_{t=0}^{T-1} \text{log}\left( \pi(\textbf{u}_t|\textbf{x}_t,\theta)\right)\right]$ as an average of $K$ samples:
\begin{align}
\nabla_\theta \mathbb{E}_{\pmb{\tau}} \approx \frac{1}{K} \sum_{k=1}^{K}  \left[ J(\pmb{\tau}^k) \nabla_\theta \sum_{t=0}^{T-1} \text{log}\left( \pi(\textbf{u}_t^k|\hat{\textbf{x}}_t^k,\theta)\right)\right]
\end{align}
where we denote the sample $k$ as a super-index. The variance of this estimation can be reduced with the aid of an action-independent baseline $b$, which does not introduce a bias \cite{Sutton}. A simple but effective baseline is the expectation of reward under the current policy, approximated by the mean over the sampled paths:
\begin{equation}\label{baseline}
    b = \hat{J}(\theta) \approx \frac{1}{K} \sum_{k=1}^KJ(\pmb{\tau}^{(k)}),
\end{equation}
which leads to:
\begin{equation}\label{Rew-Baseline}
	\nabla_\theta \hat J(\theta) \approx \frac{1}{K} \sum_{k=1}^{K}  \left[ (J(\pmb{\tau}^k)-b) \nabla_\theta \sum_{t=0}^{T-1} \text{log}\left( \pi(\textbf{u}_t^k|\hat{\textbf{x}}_t^k,\theta)\right)\right]
\end{equation}
This selection increases the log likelihood of an action by comparing it to the expected reward of the current policy. (\ref{Rew-Baseline}) is the gradient that we can now incorporate into our steepest ascent strategy. The algorithm that trains the RNN and obtains the optimal policy network is the following.

\begin{algorithm}[H]
\caption{Policy Gradient Algorithm}\label{alg:PG alg}
\smallskip

{\bf Input:} Initialize policy parameter $\theta = \theta_0$, with $\theta_0\in\Theta_0$, 
learning rate, its update rule $\alpha$, $m:=0$, the number of episodes $K$ and the number of epochs $N$.\\
{\bf Output:} policy $\pi(\cdot | \cdot ,\theta)$ and $\Theta$
\smallskip

{\bf for} m = 1,\dots, N {\bf do}
\begin{enumerate}
\item Collect $\textbf{u}_t^k , \textbf{x}_t^k$ for $T$ time steps for $K$ trajectories along with $J(\textbf{x}_T^k)$, also for $K$ trajectories.
\item Update the policy, using a policy gradient estimate $\theta_{m+1} = \theta_m + \alpha_m \frac{1}{K} \sum_{k=1}^{K}  \left[ (J(\pmb{\tau}^k)-b) \nabla_\theta \sum_{t=0}^{T-1} \text{log}\left( \pi(\textbf{u}_t^k|\hat{\textbf{x}}_t^k,\theta)\right)\right]$
\item  $m:=m+1$
\end{enumerate}
\end{algorithm}
The steps in the Algorithm \ref{alg:PG alg} are explained below.

{\bf Initialization:} The RNN policy network and its weights $\theta$ are initialized, along with the algorithm's hyperparameters such as learning rate, number of episodes and number of epochs.

{\bf Training loop:} The weights on the RNN are updated by a policy gradient scheme for a total of $N$ epochs. In {\bf Step 1} $K$ trajectories are computed, each trajectory consists of $T$ time steps, and states and control actions are collected. In {\bf Step 2} the weights of the RNN are updated based on the policy gradient framework. In {\bf Step 3}, either the algorithm terminates or returns to Step 1.

\subsection{Reinforcement Learning for Bioprocess Optimization under Uncertainty}
The methodology presented aims to overcome plant-model mismatch in uncertain dynamic systems, a usual scenario in bioprocesses. It is common to construct simple deterministic models according to a hypothesized mechanism, however the real system is more complex and presents disturbances. We propose the following methodology to address this problem (following from Algorithm \ref{alg:B2B}). 

{\bf Step 0, Initialization:} The algorithm is initialized by considering an initial policy network (e.g. RNN policy network) with untrained parameters $\theta_0$. 

{\bf Step 1, Preliminary Learning (Off-line):} It is assumed that a preliminary mechanistic model can be constructed from previous existing process data, hence, the policy learns this preliminary mechanistic model. This is done by running Algorithm~\ref{alg:PG alg} in a simulated environment by the mechanistic model. This allows the policy to incorporate previously believed knowledge about the system. The policy will therefore end with an optimal control policy for the mechanistic model. The termination criteria can be defined either by the designer or by the difference from the solution of the OCP, since the process model is known. 

Given that the experiments are in silico, a large number of episodes and trajectories can be generated that corresponds to different actions from the probability distribution of $\textbf{u}_t$, and a specific set of parameters of the RNN, respectively. The resulting control policy is a good approximation of the optimal policy. Notice that if a stochastic preliminary model exists, this approach can immediately exploit it, contrary to traditional NMPC approaches. This finishes the in silico part of the algorithm, subsequent steps would be run in the true system. Therefore, emphasis after this step is given on sampling as least as possible, as every new sample would result in a new batch run from the real plant.

{\bf Step 2-3, Transfer Learning:} The policy could directly be retrained using the true system and adapt all the weights according to the Reinforce algorithm. However, this may result in undesired effects. The control policy proposed in this work has a deep structure, as a result a large number of weights could be present. Thus, the  optimization to update the policy may easily be stuck in a low-quality local optima or completely diverge. 
To overcome this issue the concept of transfer learning is adopted. In transfer learning, a subset of training parameters is kept constant to avoid the use of a large number of epochs and episodes, applying knowledge that has been stored in a different but related problem. This technique is originated from the task of image classification, where several examples exists, $e.g.$ in \cite{Krizhevsky2012}, \cite{ILSVRC15}, \cite{Donahue2013}.

\begin{figure}[H]
    \centering
    \includegraphics[scale = 0.5]{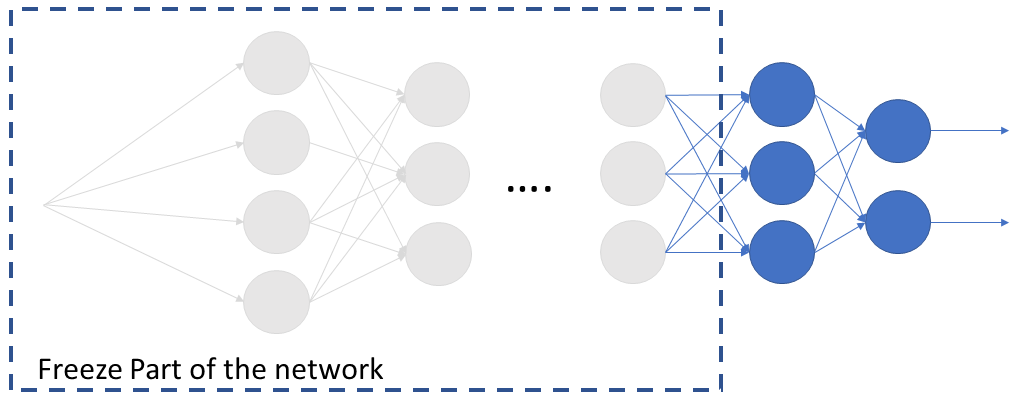}
    \caption{Part of the network is kept frozen to adapt to new situations more efficiently}
    \label{fig:rnn-frozen}
\end{figure}

Using transfer learning, the current work only retrained the last hidden layers, and the policy is able to adapt to new situation without losing previously obtained knowledge, as shown in ~Fig.\ref{fig:rnn-frozen}. Alternatively, an additional set of layers could be added on the top of the network.

{\bf Step 4, Transfer Learning Reinforce (On-line):} In this step, Algorithm~\ref{alg:PG alg} is applied again, but now, on the true system. This algorithm aims to maximize a given reward (e.g. product concentration, economic objective)

{\bf Step 5:} Terminate policy update and output $\Theta$ that defines the optimal RNN policy.

\begin{figure}[H]
    \centering
    \includegraphics[scale=0.4]{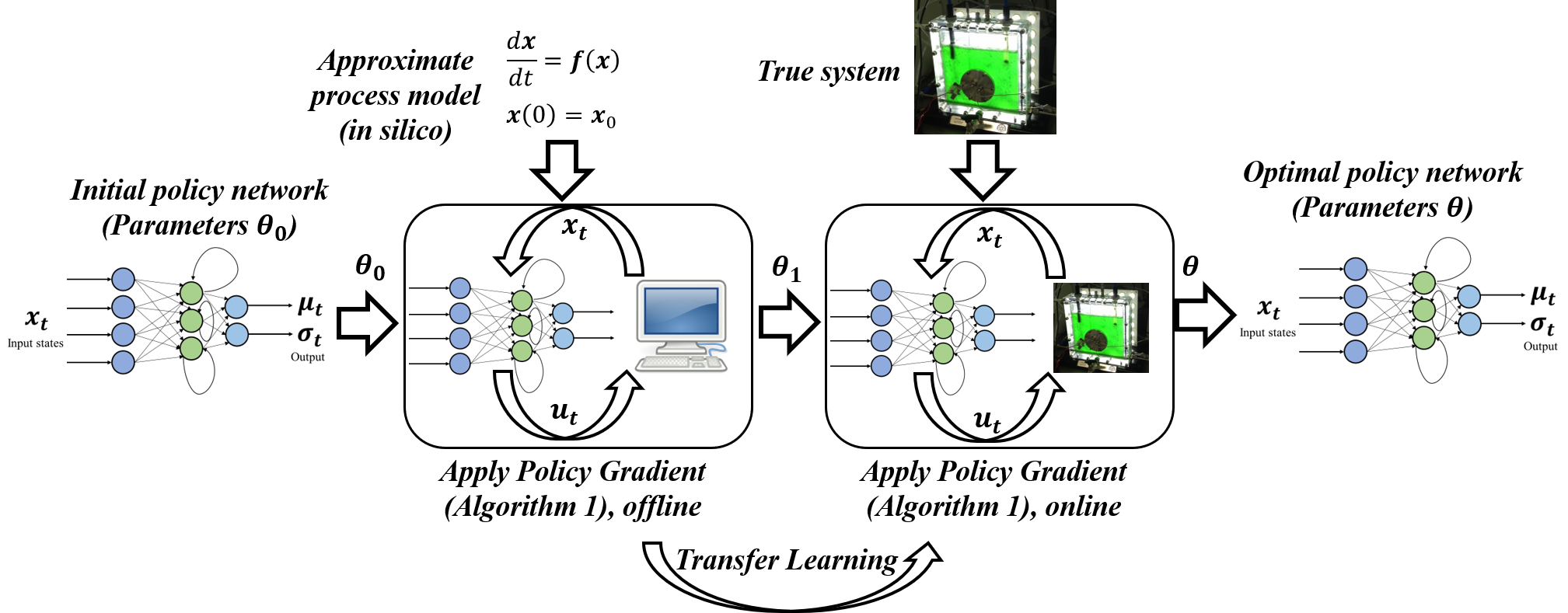}
    \caption{Batch-to-Batch algorithm (Algorithm 2)}
    \label{fig:BtB}
\end{figure}

The methodology is described in Algorithm~\ref{alg:B2B} and depicted in Fig.~\ref{fig:BtB}.
\begin{algorithm}[H]
\caption{Batch to batch algorithm}\label{alg:B2B}
\smallskip
{\bf Input:} Initialize the set of policy parameter $\Theta_0$, 
learning rate and its update rule $\alpha$, epochs $N:=N_0$, maximum number of epochs $N_{max}$, epochs for the \emph{true} system $N_{true}$, episodes $K_0$, and episodes for the \emph{true} system $K$ with $K_0>>K$.
\smallskip
\begin{enumerate}
\item \textbf{while}  
$N \leq N_{max}$ \textbf{do}:
\begin{enumerate}
\item Apply Algorithm~\ref{alg:PG alg} using an approximate model and get the trained parameters $\hat{\Theta}_0$ using $N$ epochs and $T_0$ episodes.
\item increase $N$.
\end{enumerate}
    \item $\hat{\Theta}_1$:=$\hat{\Theta}_{0}$ Initial values of the parameters $\hat{\Theta}_1$ are set as those identified in $\textbf{Step 1}$
    \item {\bf Transfer Learning:} Pick $\hat{\Theta}_1^*\subset \hat{\Theta}_1$ to be constant 
    \item For $i = 1,\dots, N_{true} $ \textbf{do}:
\begin{enumerate}
    \item Apply Algorithm~\ref{alg:PG alg} on the true system and get the trained parameters $\hat{\Theta}_1$ for one epoch and $K$ episodes.
\end{enumerate} 
\item $\Theta = \hat{\Theta}_1$
\end{enumerate}
{\bf Output:} Preliminary trained policy network with parameters $\Theta$ that takes states as inputs (e.g. $\textbf{x}_t$) and outputs the statistical parameter (e.g. $\pmb{\mu}_{t+1}$, $\pmb{\sigma}_{t+1}$) of an optimal stochastic action.

{\bf Note:} We denote $\Theta_0$ as the set of parameters of the RNN before any training, $\hat{\Theta}_0$ the set of parameters after the training in {\bf Step 1}. $\hat{\Theta}_1$ denotes the set of parameters passed along to the training by the \emph{true} system, and subsequent set of parameters during {\bf Step 4} as $\hat{\Theta}_i$, where $i$ is the current epoch.
\end{algorithm}
%
\section{Computational Case Studies}

In this section three case studies are presented to illustrate the effectiveness of the proposed batch-to-batch strategy. 
Our strategy is applied to 3 fed-batch bioreactors, where the objective is to maximize the concentration of a  target product ($y_2$ or $c_q$) at the end of the batch time, using light and an inflow rate ($u_1$ or $I$ and $u_2$ or $F_N$) as manipulated variables. 
\subsection{Case Study 1 - Ordinary Differential Equations}

In the first case study, the ``real'' photo-production system (plant) is described by the following equations plus an additional random disturbance:
\begin{align}\label{plant}
        &\dfrac{d y_1}{d t} = -({u}_1 + 0.5~{u}_1^2)  y_1 + 0.5 \dfrac{{u}_2 y_2}{(y_1 + y_2)} \\
        &\dfrac{d y_2}{d t} = {u}_1~y_1 - 0.7 {u}_2~y_1
        \end{align}
where ${u}_1$, ${u}_2$ and $y_1$, $y_2$ are the manipulated variables and the outlet concentrations of the reactant and product, respectively. The batch operation time course is normalized to 1. Additionally, a random disturbance is assumed, which is given by a Gaussian distribution with mean value 0 and standard deviation 0.02 on the states $y_1$ and $y_2$. We discretize the time horizon into 10 intervals of the dimensionless time, with one constant control input in each interval, resulting in a total of 20 control inputs.

The exact model is usually not known, and a simplified deterministic model is assumed according to some set of parameters. This preliminary model, given in (\ref{InModel1} - \ref{InModel2}),  is utilized in an extensive offline training in order to construct the control policy network. As illustrated in the previous section~\ref{REA}, there is a potential to have a close approximation of the solution of the OCP using the RNN-Reinforce.
\begin{align}
        &\dfrac{d y_1}{d t} = -({u}_1 + 0.5~{u}_1^2)
        y_1 + {{u}_2} \label{InModel1} \\
        &\dfrac{d y_2}{d t} = {u}_1~y_1 - {u}_2~y_1 \label{InModel2}
\end{align}
The training consist of 100 epochs and 800 episodes using the simplified model to search the optimal control policy that maximizes the reward of (\ref{InModel1} - \ref{InModel2}) in this case the concentration of $y_2$ at the final instance (\ref{Reward_1_s}).

\begin{equation}\label{Reward_1_s}
    \begin{split}
             R_t &= 0,~~~~t\in \{0, T-1\}\\
     R_T &= y_2(T).    
    \end{split}
\end{equation}

The control actions are constrained to be in the interval $[0,5]$. The control policy RNN is designed to contain 2 hidden layers, each of which comprises 20 neurons embedded by a hyperbolic tangent activation function. It was observed that 2 hidden layers are sufficient to approximate the optimal control policy, however there is the potential to use deeper structures with more layers for more complex systems. Furthermore, we employed two policy networks instead of one for simplicity. This approach assumes that the two manipulated variables are independent resulting in diagonal variance. 

\begin{figure}[H]
\centering
\begin{subfigure}{.5\textwidth}
  \centering
  \includegraphics[width=1.\linewidth]{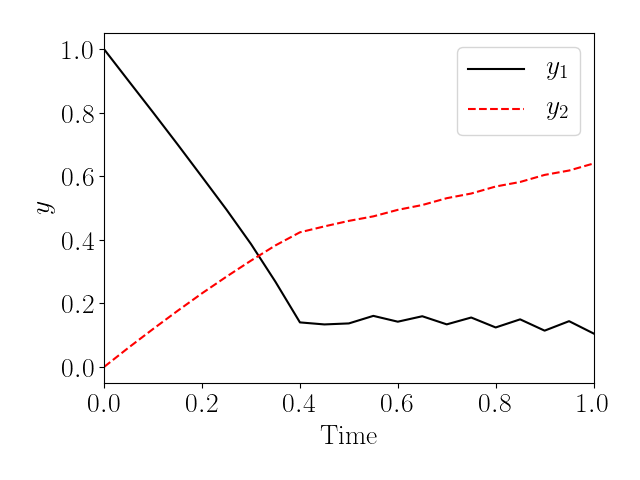}
  \caption{Solution of  $y$ for the nominal system}
  \label{fig:input_1}
\end{subfigure}%
\begin{subfigure}{.5\textwidth}
  \centering
  \includegraphics[width=1.\linewidth]{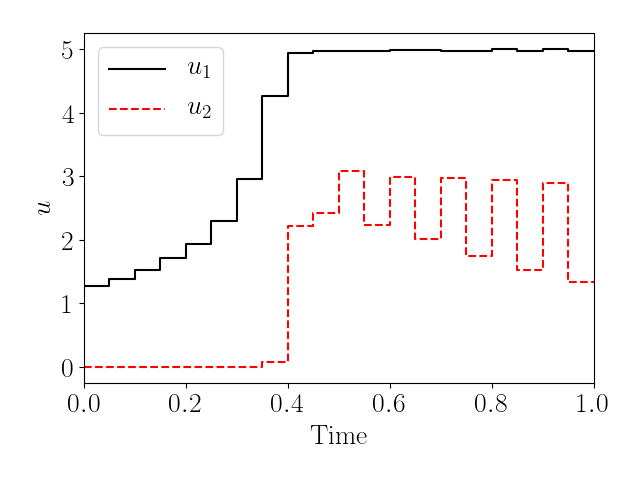}
  \caption{Solution of $u$ for the nominal system}
  \label{fig:input_2}
\end{subfigure}
\caption{The time trajectories of of the output variables of  the approximate model and the piecewise constant control actions associated with the preliminary trained policies.}
\label{fig:sol12}
\end{figure}

\begin{figure}[H]
\centering
\includegraphics[scale = 0.5]{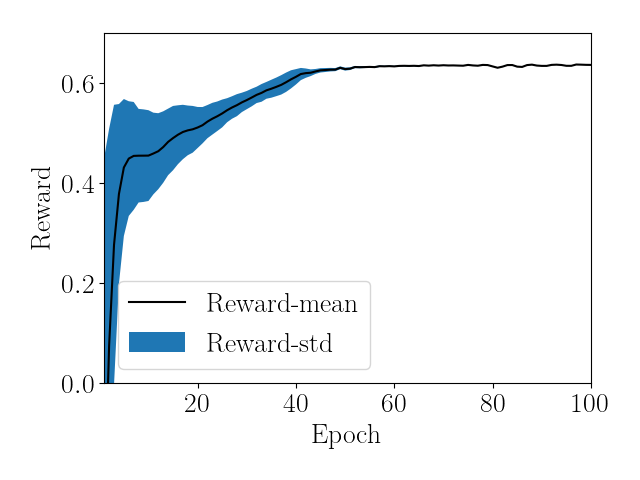}
\caption{The reward computed for the approximate model for each epoch}
\label{fig:sol22}
\end{figure}

The algorithm is implemented in Pytorch version 0.4.1.  Adam~\citep{Kingma2014} is employed to compute the network's parameter values using a step size of $10^{-2}$ with the rest of hyperparameters at their default values. 
After the training, using the simplified model the reward has the same value with the one computed by the optimal control problem, as expected. It should be noted that the computation cost of the control action using the policy is insignificant since it only requires the evaluation of the corresponding RNN, and does not depend directly on the complexity or the number of variables. In contrast, the solution of the OCP scale very badly with respect to both the complexity and the number of variables.
Precisely, the maximum rewards for RL and OCP for both cases is 0.64. The reward for its epoch is illustrated in Fig.~\ref{fig:sol22} and the process trajectories after the final update of the policy networks are shown in Fig.~\ref{fig:sol12}.

\begin{figure}[H]
\centering
\begin{subfigure}{.5\textwidth}
  \centering
  \includegraphics[width=1.\linewidth]{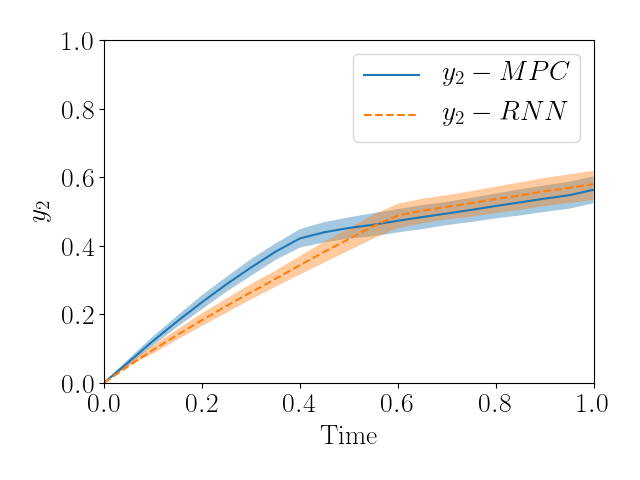}
  \caption{Solution of  $y_2$ for the ``real'' system}
  \label{fig:input_1}
\end{subfigure}%
\begin{subfigure}{.5\textwidth}
  \centering
  \includegraphics[width=1.\linewidth]{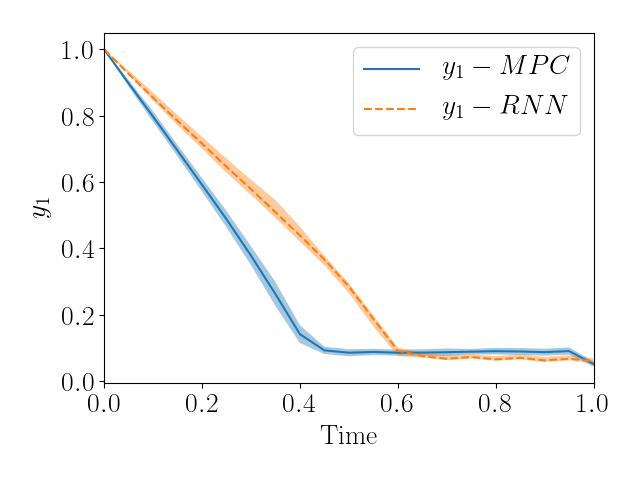}
  \caption{Solution of $y_1$ for the ``real'' system}
  \label{fig:input_2}
\end{subfigure}
\caption{The time trajectories produced by the real plant using our approach (dash) and NMPC (solid).}
\label{fig:sol11}
\end{figure}

\begin{figure}[H]
\centering
\includegraphics[scale = 0.6]{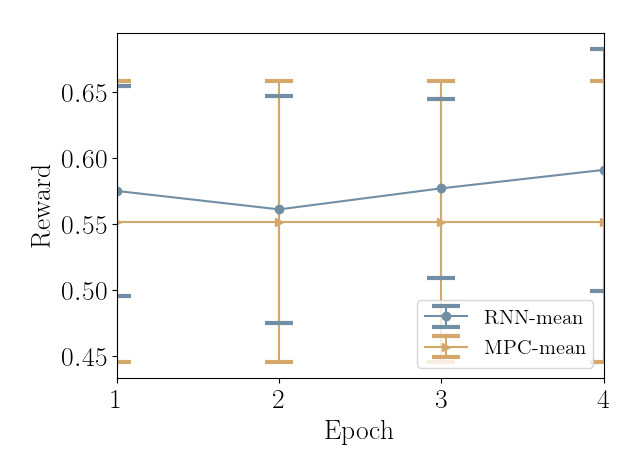}
\caption{The reward computed by the updated training using the plant (``real" system) for each epoch (circle) and the average performance of the NMPC (triangle) with 2 times the standard deviation.}
\label{fig:sol21}
\end{figure}

Fig.~\ref{fig:sol22} shows that the reward has a large variance at the beginning of the training but is undetectable at the end. This can be explained as the trade-off between exploration and exploitation, where initially there is a lack of information and policy explores possible control actions, while at the end the policy exploits the available information. 
This policy can be considered as an initialization of the Reinforce algorithm which uses transfer learning to incorporate new knowledge gained from the true plant (Steps 3-5 in Algorithm \ref{alg:B2B}). New data-sets from 25 batches are used ({\it i.e.} 25 real plant epochs) to update the true plant's RL policy. The solution after only 4 epochs is 0.591 while the stochastic-free optimal solution identified using the {\it unknown} (complex) model of the plant is 0.583. This results show that the stochastic nature of the system can also affect the performance.  
The reward for each epoch is depicted in Fig.~\ref{fig:sol21} and the process trajectories after the last epoch are depicted in Fig.~\ref{fig:sol11}. Notice that even before having any interaction with the ``real'' system the proposed approach has a superior performance than NMPC. This is because RL directly accounts for the effect of future uncertainty and its feedback in a proper closed-loop manner, whereas NMPC assumes open-loop control actions at future time points in the prediction.

There is a variation on the results after the last batch upon which the policy is updated. This makes sense, since the system appears to have some additive noise (i.e. Gaussian disturbance) and the policy maintains its stochastic nature. 

\begin{figure}[H]
\centering
\begin{subfigure}{.5\textwidth}
  \centering
  \includegraphics[width=1.\linewidth]{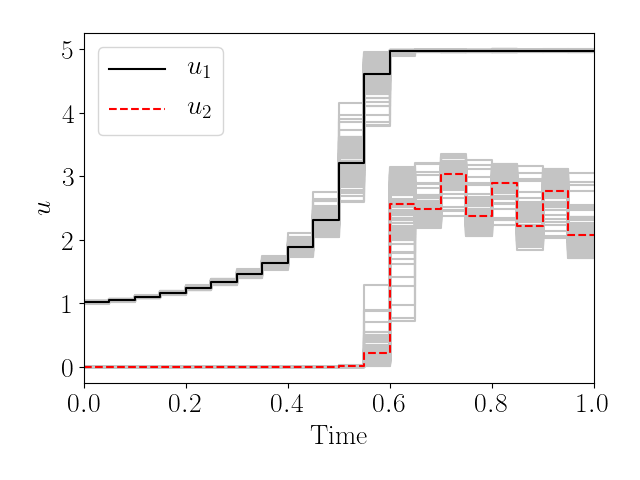}
  \caption{Solution of RNN for the ``real'' system}
  \label{fig:INP1}
\end{subfigure}%
\begin{subfigure}{.5\textwidth}
  \centering
  \includegraphics[width=1.\linewidth]{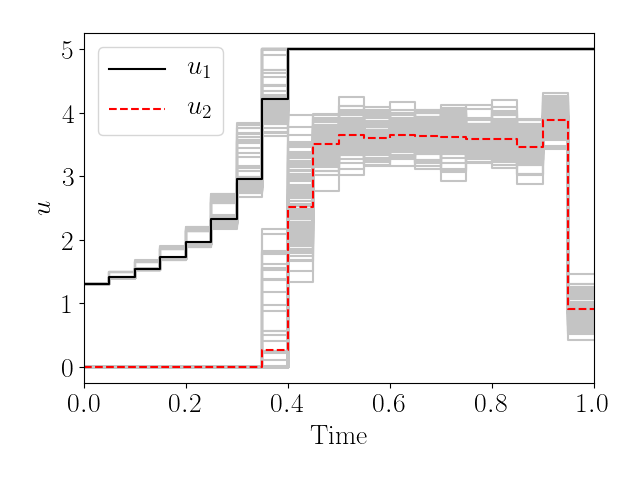}
  \caption{Solution of NMPC for the ``real'' system}
  \label{fig:INP2}
\end{subfigure}
\caption{Comparison of the time trajectories of the piecewise constant control actions between our approach (left) and NMPC (right)}
\label{fig:inputs_1}
\end{figure}

The results are also compared with the use of NMPC using shrinking horizon. The results can be seen in Fig.~\ref{fig:sol21}, where 100 Monte-Carlo simulations were conducted. The optimization using our approach appears to be superior to the one given by the NMPC, showing the significance of our result.  Furthermore, it should be noted that the performance of our proposed policy is better even in epoch 1, before the adaptation is started. In addition, in Fig.~\ref{fig:inputs_1}, the comparison between the control inputs of that are computed using our approach and the NMPC. 

\subsection{Case Study 2 - Stochastic Differential Equations}
In this case study the same type of reaction is assumed to follow a stochastic differential equations: 
    \begin{align*}
        &d y_1 = \left[-({u}_1 + 0.5~{u}_1^2)  y_1 + 0.5 \dfrac{{u}_2 y_2}{(y_1 + y_2)}\right]dt \\
        &{d y_2} = \left[{u}_1~y_1 - 0.7 {u}_2~y_1\right]dt + \left[0.1 \sqrt{y_1}\right]dW
    \end{align*}
where $W$ is Wiener stochastic process. The simplified model is assumed to be the same with the previous case study. As a result the same policy that is trained off-line is used here. The purpose of this case study is to observe how the same policy can adapt in different environments. Now the model that describes the real system is not only structurally different, but also stochastic in nature.

\begin{figure}[H]
\centering
\begin{subfigure}{.5\textwidth}
  \centering
  \includegraphics[width=1.\linewidth]{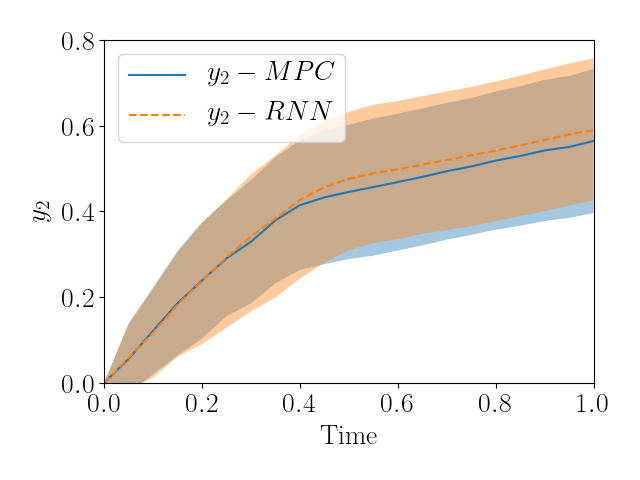}
  \caption{Solution of  $y_2$ for the ``real'' system}
  \label{fig:input_1}
\end{subfigure}%
\begin{subfigure}{.5\textwidth}
  \centering
  \includegraphics[width=1.\linewidth]{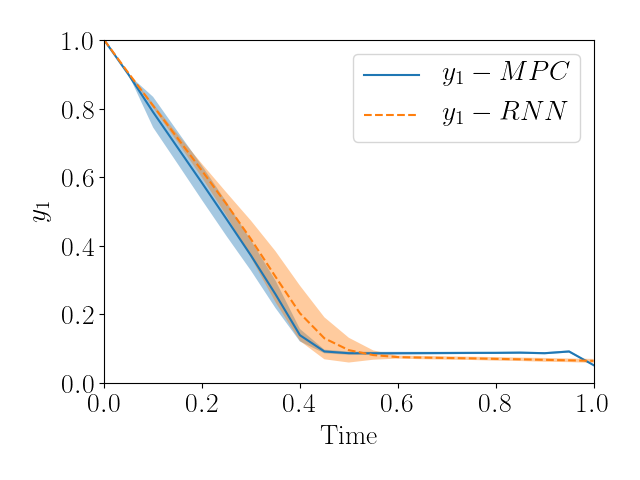}
  \caption{Solution of $y_1$ for the ``real'' system}
  \label{fig:input_2}
\end{subfigure}
\caption{The time trajectories produced by the real plant using our approach (dash) and NMPC (solid).}
\label{fig:SDE_Result}
\end{figure}

\begin{figure}[H]
\centering
\begin{subfigure}{.5\textwidth}
  \centering
  \includegraphics[width=1.\linewidth]{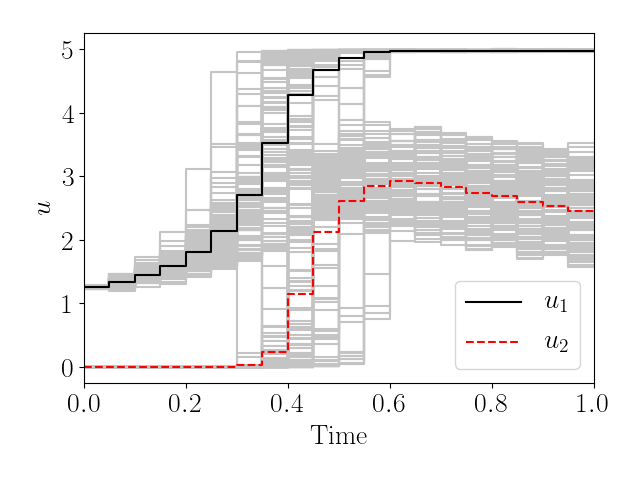}
  \caption{Solution of RNN for the ``real'' system}
  \label{fig:inut_1}
\end{subfigure}%
\begin{subfigure}{.5\textwidth}
  \centering
  \includegraphics[width=1.\linewidth]{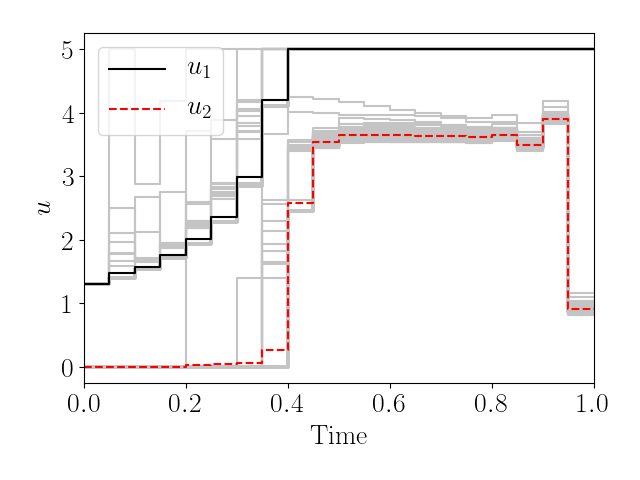}
  \caption{Solution of NMPC for the ``real'' system}
  \label{fig:inut_2}
\end{subfigure}
\caption{Comparison of the time trajectories of the piecewise constant control actions between our approach and NMPC}
\label{fig:inputs_2}
\end{figure}

The same hyperparmeters and networks are utilized for the policies in both stages, in order to show that the same policy can adapt to different environments successfully. 
\begin{figure}[H]
    \centering
    \includegraphics[scale = 0.6]{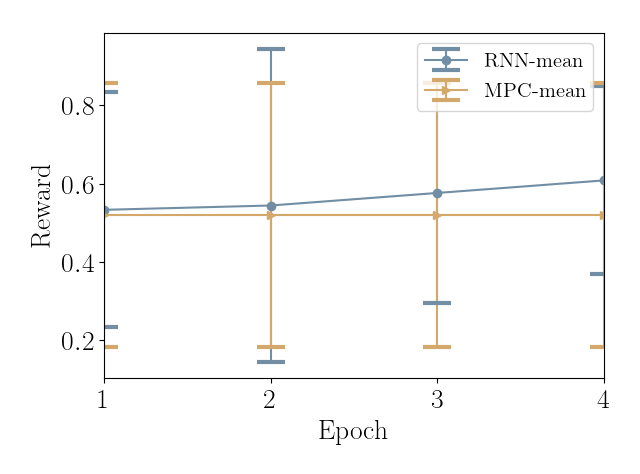}
    \caption{The reward computed by the updated training using the plant (``real" system) for each epoch (circle) and the average performance of NMPC (triangle)  with 2 times the standard deviation.}
    \label{fig:SDE_Result1}
\end{figure}
The same validation is conducted here using 100 Monte-Carlo simulations. Through comparison, our approach is  found to be superior to the NMPC.  In this case, our proposed algorithm adapts more rapidly to the new conditions, reducing significantly the requirement for a large number of episodes and epochs, as it can be seen in Fig.~\ref{fig:SDE_Result1}. This is attributed to the systematic transfer learning proposed in our algorithm. The computationally intensive part has been shifted off-line where the preliminary inaccurate model was used to train the policy. Then the (deep) recurrent neural network adapts successfully to the new environment that consists of a system of stochastic differential equations. The comparison is also depicted in Fig.~\ref{fig:SDE_Result}. This result is also observed in the previous case study where the stochastic part of the physical system has a  different nature. The control inputs are depicted in Fig.~\ref{fig:inputs_2}.

It should be noted that in both case studies the NMPC produced very similar control actions, with the only difference being the variance, compared to our approach which shapes the control actions to fit the needs of the different dynamics and uncertainty. 

The methods used in the Reinforce algorithm usually require substantial number of  episodes and epochs, therefore a good initial solution in combination with transfer learning is paramount so that Step 4 can be completed with a few batch-to-batch runs. In order to keep the problem realistic, only a small number of batches is utilized in Step 2-3 to refine the policy network.

\subsection{Case Study 3 - Nonsmooth Model}

The last case study in this paper focuses on the photo-production of phycocyanin synthesized by cyanobacterium
Arthrospira platensis. Phycocyanin is a high-value bioproduct and its biological function is to enhance the photosynthetic
efficiency of cyanobacteria and red algae. It has applications as a natural
colorant to replace other toxic synthetic pigments in both food and cosmetic production. Additionally, the pharmaceutical industry considers it as beneficial because of its unique antioxidant, neuroprotective, and anti-inflammatory properties.

Both the ``real'' and simplified model in this case are considered to be nonmooth. Due to different growth phases, nonmooth behaviour is observed for the physical system. To accommodate this difficulty, switching functions have been proposed~\citep{Antonio-nonsmooth, Dongda-nonsmooth}. In this work the nonsmooth behaviour is modelled using a $\sign(\cdot)$ function. 
The ``real'' dynamic system consists of three nonsmooth ODEs describing the evolution of the concentration of biomass ($X$), nitrate ($N$), and product ($q$). The dynamic model is based on Monod kinetics, which describes microorganism growth in nutrient sufficient cultures, where intracellular nutrient concentration is kept constant because of the rapid replenishment. We assume a fixed volume fed-batch. The manipulated variables as in the previous examples are the light intensity ($I$) and inflow rate ($F_N$). The mass balance equations are
 
 \begin{align}
     \dfrac{dc_x}{dt} &=  u_m \dfrac{I}{I + k_s + I^2/k_i}c_x \dfrac{c_N}{c_N + K_N} - u_dc_X\label{biomas}\\
     \dfrac{dc_N}{dt} &= -Y_{N/X} u_m \dfrac{I}{I + k_s + I^2/k_i}c_x \dfrac{c_N}{c_N + K_N} + F_N \label{Nitrogen}\\
     \dfrac{dc_q}{dt} &=  \begin{cases}
      k_m \dfrac{I}{I + k_{sq} + I^2/k_{iq}}c_x \dfrac{c_N}{c_N + K_N} - k_d\dfrac{c_{q}}{C_N + K_Nq} &,\text{if}\ c_N \leq 500mgL^{-1} \& c_X\geq 10 gL^{-1}\\
      0 &,\text{otherwise} \label{product},
    \end{cases}
\end{align}
 
 where the parameters are given in Table~\ref{tab:paramvalues}. The real physical system consists of additive disturbance
 
 \begin{align}
          w(t) &= sin(t) \sigma_d + \sigma_n\\
          \sigma_d &= diag(4\times 10^{-3}, 1., 1\times10^{-7})\\
          \sigma_n &\sim \mathcal{N}(\textbf{0}, \sigma_d),
 \end{align}

and measurement noise

  \begin{align}
          noise(t) \sim \mathcal{N}(\textbf{0}, diag(4\times 10^{-4}, .1, 1\times10^{-8})).
 \end{align}
 
Additionally, uncertainty is assumed for the initial concentration, where

\begin{align}
    \begin{bmatrix}
    c_x(0) & c_N(0) & c_q(0)
    \end{bmatrix} \sim \mathcal{N}(    \begin{bmatrix}
    1. & 150. & 0.
    \end{bmatrix}, diag(1\times 10^{-3}, 22.5, 0.)).
\end{align}
 
 The reward is additionally penalized by the change of the control actions $\textbf{u}(t) = \left[ I, F_N\right]^{T}$. As a result the reward is given as:
 \begin{equation}\label{reward_bio}
 \begin{split}
     R_t &= -\Delta \textbf{u}_t^T diag(3.125\times10^{-8}, 3.125 \times 10^{-6})  \Delta \textbf{u}_t^T,~~~~~~~ t\in \{0, T-1\}\\
     R_T &= c_q(T),     
 \end{split}
 \end{equation}
 where $\Delta \textbf{u}_t = \textbf{u}_t - \textbf{u}_{t-1}$.
 
The simplified deterministic model is assumed without the noise or the additive disturbance. This preliminary model, is utilized in an extensive offline training in order to construct the control policy network.  As illustrated in the previous section 3.4, there is a potential to have a close approximation of the solution of the OCP using RNN-Reinforce.

The training consists on 100 epochs and 500 episodes and the optimal control policy that maximizes the reward in equation~\ref{reward_bio}. The control actions are constrained to be in the interval $0\leq F_N \leq 40$ and $120 \leq T\leq 400$. The control policy RNN is designed to contain 4 hidden layers, each of which comprises 20 neurons embedded by a leaky rectified linear unit (ReLU)  activation function. Furthermore, in this case a unified policy network with diagonal variance is utilized such that the control actions share memory and the previous states are used from the RNN (together the current measured states).

\begin{figure}[H]
\centering
\includegraphics[scale = 0.6]{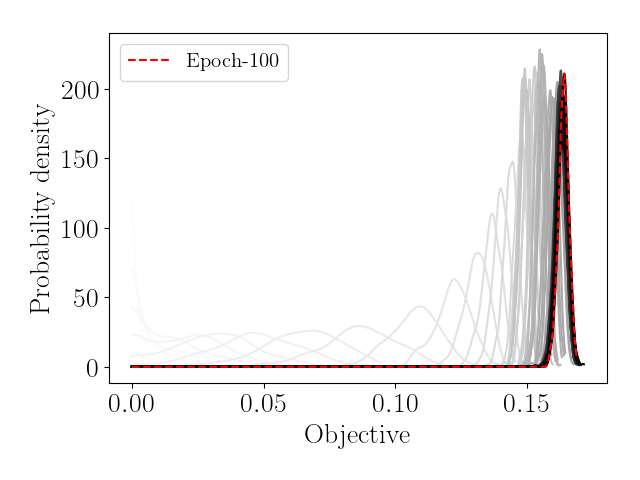}
\caption{The reward computed for the approximate model for each epoch}
\label{fig:epoch-bio}
\end{figure}

The algorithm is implemented in Pytorch with the same configurations. It  should  be  noted  that  the  computational cost of computing the control action online is insignificant since it only requires the evaluation of the corresponding RNN, and does not depend directly on the complexity or the number of variables. In contrast, the solution of the OCP scales badly in this case due to the presence of integer variables. The reward for each epoch is depicted in Fig~\ref{fig:epoch-bio}. In this case the probability density is shown due to the nonsmoothness of the model, that may result in multiple peaks. The lines  are faded out towards earlier epochs. Additionally, there is no guarantee of global optimality in the current work, as a result the reward may get stuck in other local minima. It should be noted that in this case the uncertain initial conditions are applied during this training phase.

\begin{figure}[H]
\centering
\begin{subfigure}{.5\textwidth}
  \centering
  \includegraphics[width=1.\linewidth]{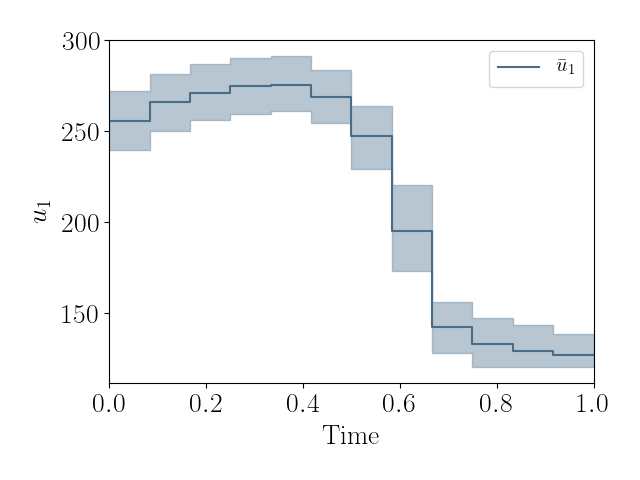}
  \caption{Nominal solution for $I$ }
  \label{fig:u1_nominal}
\end{subfigure}%
\begin{subfigure}{.5\textwidth}
  \centering
  \includegraphics[width=1.\linewidth]{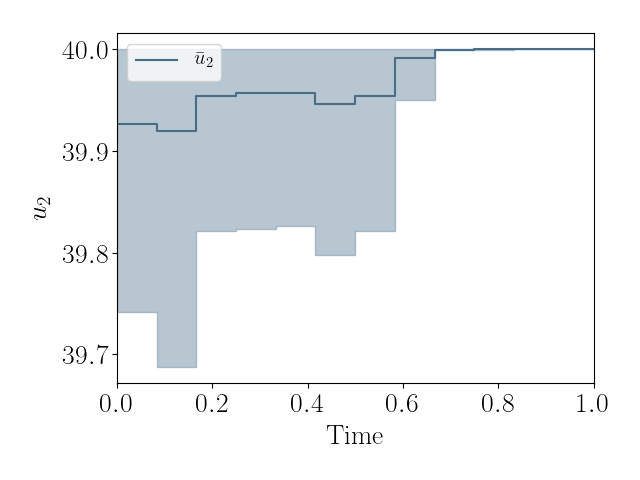}
  \caption{Nominal solution for $F_N$}
  \label{fig:u2_nominal}
\end{subfigure}
\caption{Solution for control actions of the nominal system using RNN}
\label{fig:unominal_bio}
\end{figure}

\begin{figure}
\begin{subfigure}{.5\linewidth}
\centering
\includegraphics[width=.9\linewidth]{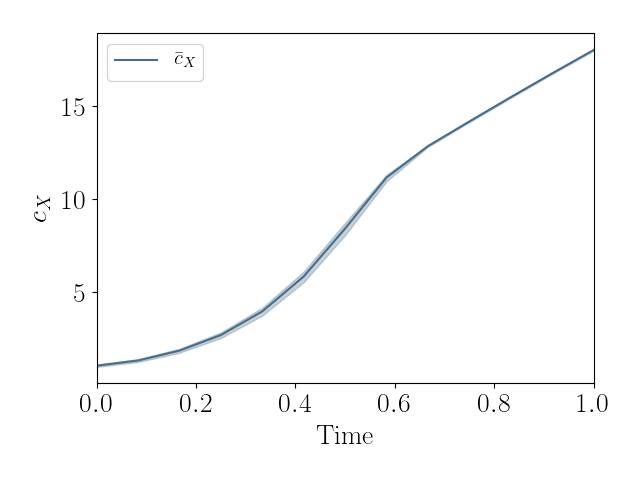}
\caption{The nominal solution for $c_X$}
\label{fig:cx}
\end{subfigure}%
\begin{subfigure}{.5\linewidth}
\centering
\includegraphics[width=.9\linewidth]{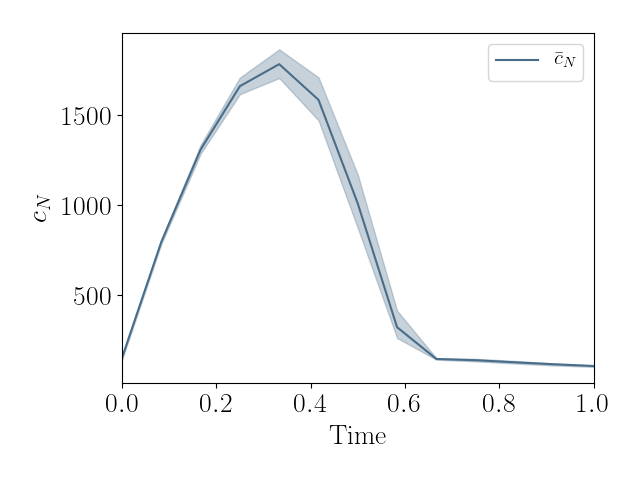}
\caption{The nominal solution for $c_N$}
\label{fig:cn}
\end{subfigure}\\[1ex]
\begin{subfigure}{\linewidth}
\centering
\includegraphics[width=.5\linewidth]{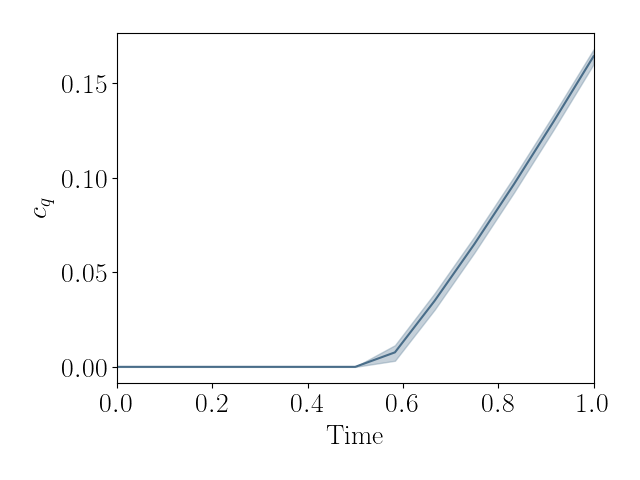}
\caption{The nominal solution for $c_q$}
\label{fig:cq}
\end{subfigure}
\caption{Comparison of the responses when RNN and NMPC are applied to the ``real'' physical system.}
\label{fig:total_nominal}
\end{figure}

 The nominal control actions are depicted in Fig.~\ref{fig:unominal_bio}, where the shaded areas are the 98\% and 2\% percentiles.  The corresponding states are depicted in Fig.\ref{fig:total_nominal} with their 98\% and 2\% percentiles. The nominal behaviour is subject to the corresponding initial conditions since no other uncertainty is taken into account in the offline procedure. The probability density of the product $c_q$ has clearly only one peak, this is shown in Fig.~\ref{fig:cq}.

\begin{table}[htbp]
  \centering
  \caption{Parameter values for physical system (\ref{biomas} - \ref{product})}
    \begin{tabular}{lrl}
        \hline
    \multicolumn{3}{c}{Parameter Values} \\
        \hline
    $u_m$  & 0.0572 &  $h^{-1}$ \\
    $u_d$  & 0.0     &  $h^{-1}$ \\
    $K_N$  & 393.1 & $mg/L$ \\
    $Y_{NX}$ & 504.1 & $mg/g$ \\
    $k_m$  & 0.00016 & $mg/g/h$ \\
    $k_d$  & 0.281 & $h^{-1}$ \\
    $k_s$  & 178.9 & $\mu mol/m^{2}/s$ \\
    $k_i$  & 447.1 & $\mu mol/m^{2}/s$ \\
    $k_{sq}$ & 23.51 & $\mu mol/m^{2}/s$ \\
    $k_{iq}$ & 800   & $\mu mol/m^{2}/s$ \\
    $K_{NP}$ & 16.89 & $mg/L$\\
            \hline
    \end{tabular}%
  \label{tab:paramvalues}%
\end{table}%

As in the previous case studies, the results are compared with the use of NMPC using shrinking horizon. The optimization is a mixed integer nonlinear programming problem (MINLP). Local optimization is used in order to be numerically tractable. Orthogonal collocation is implemented and integer variables have been used to model the switches. It should be noted that this MINLP takes 2-4 mins to be solved. 

\begin{figure}
\begin{subfigure}{.5\linewidth}
\centering
\includegraphics[width=.9\linewidth]{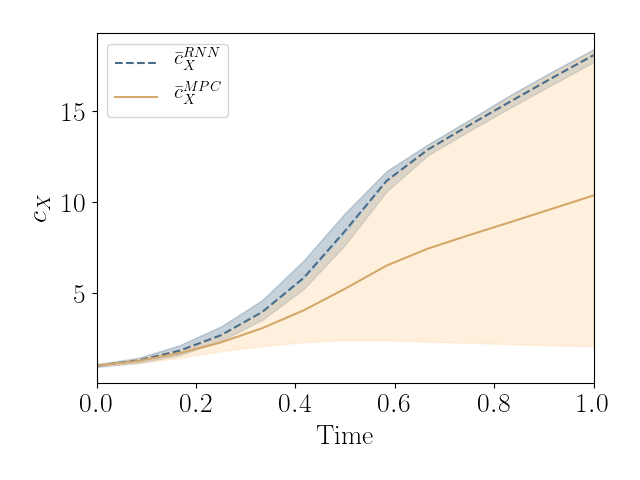}
\caption{The solution for $c_X$}
\label{fig:sub1}
\end{subfigure}%
\begin{subfigure}{.5\linewidth}
\centering
\includegraphics[width=.9\linewidth]{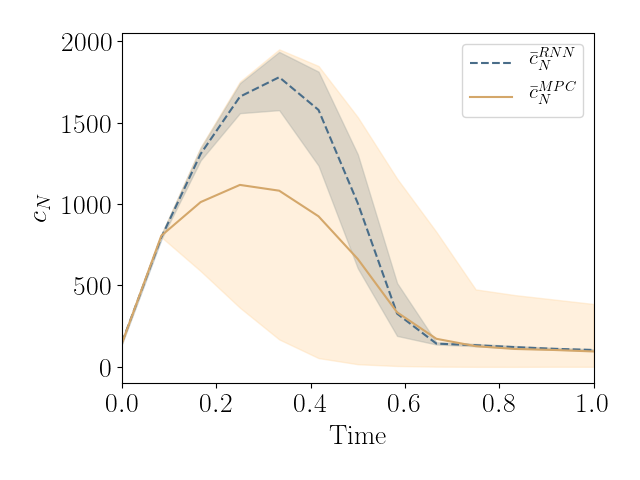}
\caption{The solution for $c_N$}
\label{fig:sub2}
\end{subfigure}\\[1ex]
\begin{subfigure}{\linewidth}
\centering
\includegraphics[width=.5\linewidth]{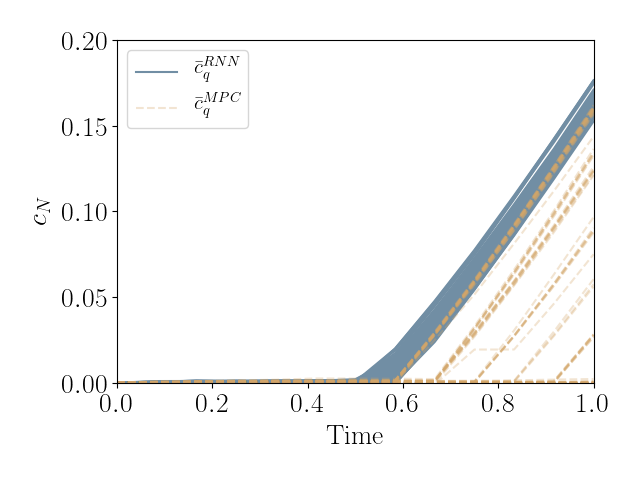}
\caption{The solution for $c_q$}
\label{fig:sub3}
\end{subfigure}
\caption{Comparison of the responses when RNN and NMPC are applied to the ``real'' physical system.}
\label{fig:states_bio}
\end{figure}

\begin{figure}[H]
\centering
\includegraphics[scale = 0.5]{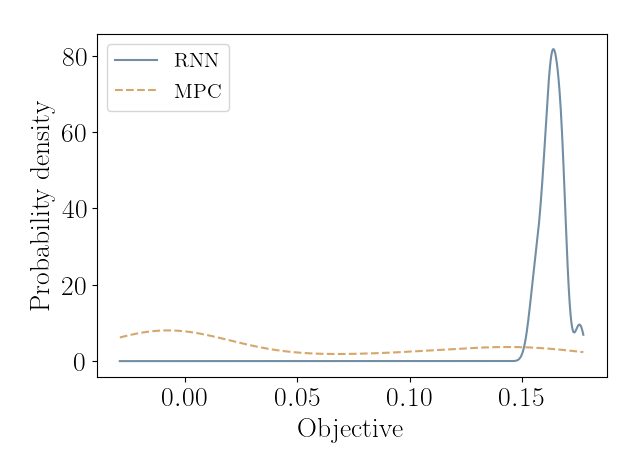}
\caption{Probability density function for the reward computed by the updated training using the plant (``real" system) for each epoch and performance of NMPC.}
\label{fig:rew-obj-bio}
\end{figure}

\begin{figure}[H]
\centering
\begin{subfigure}{.5\textwidth}
  \centering
  \includegraphics[width=.9\linewidth]{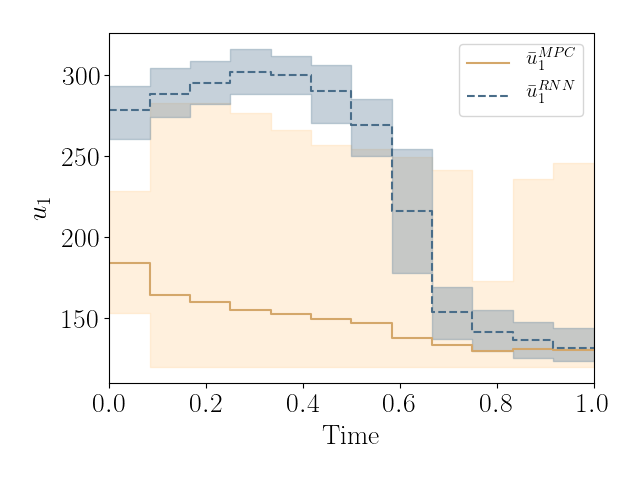}
  \caption{Solution for $I$ of the ``real'' system}
  \label{fig:input_bio1}
\end{subfigure}%
\begin{subfigure}{.5\textwidth}
  \centering
  \includegraphics[width=.9\linewidth]{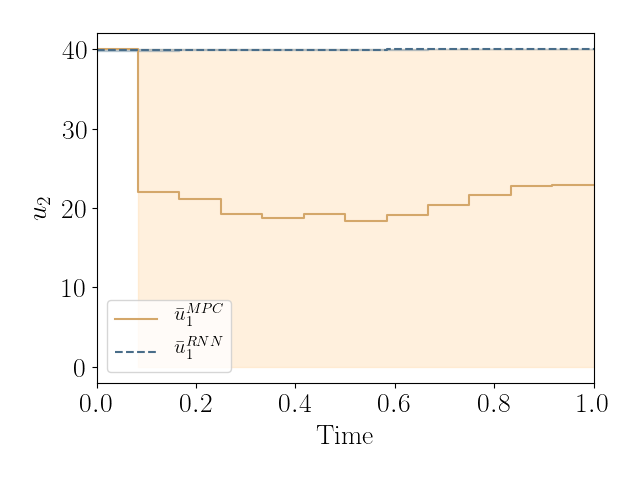}
  \caption{Solution for $F_N$ of the ``real'' system}
  \label{fig:input_bio2}
\end{subfigure}
\caption{Comparison of the time trajectories of the piecewise constant control actions between our approach and NMPC}
\label{fig:controls-real}
\end{figure}

The results can be seen in Fig.~\ref{fig:states_bio}  with  their  98\%  \&  2\%  percentile  respectively, where 100 Monte-Carlo simulations were conducted. The optimization using our approach appears to be superior to the one given by the NMPC, showing the significance of our result. After the adaptation the probability densities are depicted in Fig.~\ref{fig:rew-obj-bio}. Next, the control inputs are depicted in Fig.~\ref{fig:controls-real} with their 98\% \& 2\% percentile respectively. It is clear that the NMPC control actions have large variance compare to the ones produce by our proposed methodology. This is due to the nonsmoothness of the model and the uncertainty which the NMPC struggles with.

In addition, in Fig.~\ref{fig:controls-real}, the comparison between the control inputs of our approach and the NMPC is presented. 

\section{Conclusions and Future Work}
In this work we propose a new methodology for batch-to-batch learning by adapting Reinforcement learning techniques to uncertain and complex bioprocesses. The results reveal that it is possible to obtain a near optimal policy for a stochastic system when the true dynamics are unknown. In real systems with the absence of a true model, it is impossible to generate highly accurate datasets to train the policy network. As a result we propose a 2-stage framework where first an approximate (possibly stochastic) model is used to train the policy network. Subsequently, this policy is implemented into the \emph{true} system. In this way, there is no need for a large number of evaluations of the true system which can be costly and time consuming.

A systematic adaptation to the new environment is achieved using transfer learning. In Step 4: Transfer Learning Reinforce the policy is trained using $T<<T_0$ episodes conducting the Steps 1-3 of Algorithm~\ref{alg:PG alg}. The proposed algorithm is validated using two case studies for different nature of stochastic processes. Our proposed methodology results in a policy that overcomes the performance of the NMPC, where only simple policy evaluations are needed.  

The off-line CPU time is 3 hours, however the online implementation of the needs only 0.002 secs. This means that all the computational complexity is shifted offline and an efficient optimal control policy is constructed. One should also keep in mind that after the off-line training the solution to a nonlinear stochastic dynamical system is provided, in the form of a stochastic policy. This is a more complete and efficient solution as it is a closed-loop solution, rather than an open-loop optimization. Furthermore, a nonsmooth system was integrated with Casadi \cite{Andersson2019}, which is more time consuming than integrating a smooth dynamic system.

For both the case studies 4 epochs and 25 batches were implemented. In this work, the training was stopped after the $4^{th}$ epoch, but the training could have been continued or stopped earlier. Here, the total number of batches is $4 \times 25=100$; however, the policy for all case studies performs better from the beginning of the online implementation. This means that a smaller number of batches can be used and still outperform NMPC.

We emphasize that our considered systems contain both stochasticity and plant-model mismatch, and there is no process structure available. The optimisation of such systems is generally known to be intractable. Given the early stage of this research, there are still disadvantages of this method which must be accommodated in the future, including the robust satisfaction of constraints. In addition, there is a wide discussion regarding the safety in reinforcement learning ~\cite{DBLP:journals/corr/abs-1812-05506}, which is also a result of the difficulty of robust satisfaction of constraints.  Future work will focus on the robust satisfaction of constraints in RL methods.

The codes are available at:  \href{https://gitlab.com/Panos108/rl-with-nonsmooth}{https://gitlab.com/Panos108/rl-with-nonsmooth}

\section*{Acknowledgements}
This project has received funding from the EPSRC project (EP/P016650/1).

%
\section*{Bibliography}

\bibliography{example}
\end{document}